\documentclass{article}

\usepackage{arxiv}
\usepackage{graphicx, color}
\usepackage{amsmath,amssymb}
\usepackage[utf8]{inputenc} 
\usepackage[T1]{fontenc}    
\usepackage{hyperref}       
\usepackage{url}            
\usepackage{booktabs}       
\usepackage{amsfonts}       
\usepackage{microtype}      
\usepackage{lipsum}

\title{Lemniscate of Leaf Function}

\author{
  Kazunori Shinohara\thanks{10-3 Takiharu-cho, Minami-ku, Nagoya 457-8530, Japan} \\
  Department of Mechanical Systems Engineering\\
  Daido University\\
  10-3 Takiharu-cho, Minami-ku, Nagoya 457-8530, Japan \\
  \texttt{shinohara@06.alumni.u-tokyo.ac.jp} \\
}

\begin{document}
\maketitle

\begin{abstract}
A lemniscate is a curve defined by two foci, $F_1$ and $F_2$. If the distance between the focal points of $F_1 - F_2$ is 2a (a: constant), then any point $P$ on the lemniscate curve satisfy the equation $PF_1 \cdot PF_2 = a^2$. Jacob Bernoulli first described the lemniscate in 1694. The Fagnano discovered the double angle formula of the lemniscate(1718). The Euler extended the Fagnano's formula to a more general addition theorem(1751). The lemniscate function was subsequently proposed by Gauss around the year 1800. These insights are summarized by Jacobi as the theory of elliptic functions. A leaf function is an extended lemniscate function. Some formulas of leaf functions have been presented in previous papers; these included the addition theorem of this function and its application to nonlinear equations. In this paper, the geometrical properties of leaf functions at $n = 2$ and the geometric relation between the angle $\theta$ and lemniscate arc length $l$ are presented using the lemniscate curve. The relationship between the leaf functions $\mathrm{sleaf}_2(l)$ and $\mathrm{cleaf}_2(l)$ is derived using the geometrical properties of the lemniscate, similarity of triangles, and the Pythagorean theorem. In the literature, the relation equation for $\mathrm{sleaf}_2(l)$ and $\mathrm{cleaf}_2(l)$ (or the lemniscate functions, $\mathrm{sl}(l)$ and $\mathrm{cl}(l)$) has been derived analytically; however, it is not derived geometrically.
\end{abstract}


\section{Introduction}
\subsection{Motivation}
An ordinary differential equation (ODE) comprises a function raised to the $2n-1$ power and the second derivative of this function. Further, the initial conditions of the ODE are defined.

\begin{equation}
\frac{\mathrm{d}^2r(l) }{\mathrm{d}l^2}=-nr(l)^{2n-1} \label{1.1.1}
\end{equation}
\begin{equation}
r(0)=1 \label{1.1.2}
\end{equation}
\begin{equation}
\frac{\mathrm{d}r(0)}{\mathrm{d}l}=0\label{1.1.3}
\end{equation}

Another ODE and its initial conditions are given below.

\begin{equation}
\frac{\mathrm{d}^2 \overline{r}(\overline{l}) }{\mathrm{d}\overline{l}^2}=-n \overline{r}(\overline{l})^{2n-1} \label{1.1.4}
\end{equation}
\begin{equation}
\overline{r}(0)=0 \label{1.1.5}
\end{equation}

\begin{equation}
\frac{\mathrm{d}\overline{r}(0)}{\mathrm{d}\overline{l}}=1\label{1.1.6}
\end{equation}

The ODE comprises a function $r(l)$(or $\overline{r}(\overline{l})$) of one independent variable $l$(or $\overline{l}$) and the derivatives of this function. The variable $n$ represents a natural number ($n$ = 1, 2, 3, $\cdots$). The above equation and the initial conditions constitute a very simple ODE. However, when this differential equation is numerically analyzed, mysterious waves are generated for all natural numbers. These mysterious waves are regular waves with some periodicity and amplitude.
If these waves can be explained, they have the potential to solve various problems of nonlinear ODEs.


\subsection{Theory of Leaf Functions}
No elementary functions satisfy Eqs. (\ref{1.1.1})--(\ref{1.1.3}). Therefore, in this paper, the function that satisfies Eqs. (\ref{1.1.1})--(\ref{1.1.3}) is defined as cleaf$_n(l)$. Function $r(l)$ is abbreviated as $r$. By multiplying the derivative $dr/dl$ with respect to Eq. (\ref{1.1.1}), the following equation is obtained.

\begin{equation}
\frac{\mathrm{d}^2r }{\mathrm{d}l^2} \frac{\mathrm{d}r }{\mathrm{d}l}=-nr^{2n-1}\frac{\mathrm{d}r }{\mathrm{d}l} \label{1.2.1}
\end{equation}

The following equation is obtained by integrating both sides of Eq. (\ref{1.2.1}).

\begin{equation}
\frac{1}{2} \left( \frac{\mathrm{d}r }{\mathrm{d}l} \right)^2=-\frac{ 1 }{ 2 } r^{2n} +C\label{1.2.2}
\end{equation}

Using the initial conditions in Eqs. (\ref{1.1.2}) and (\ref{1.1.3}), the constant $C = \frac{1}{2}$ is determined. The following equation is obtained by solving the derivative $dr/dl$ in Eq. (\ref{1.2.2}).

\begin{equation}
\frac{\mathrm{d}r }{\mathrm{d}l} = \pm \sqrt{1- r^{2n}} \label{1.2.3}
\end{equation}

We can create a graph with the horizontal axis as the variable $l$ and the vertical axis as the function $r$. Because function $r$ is a wave with a period, the gradient $dr/dl$ has positive and negative values, and it depends on domain $l$. In the domain, $0 \leqq l \leqq \frac{1}{2} \pi_n$(See Appendix A for the constant $\pi_n$), the above gradient $dr/dl$ becomes negative.

\begin{equation}
{\mathrm{d}l} = - \frac{\mathrm{d}r }{ \sqrt{1- r^{2n}} } \label{1.2.4}
\end{equation}

The following equation is obtained by integrating the above equation from 1 to $r$.

\begin{equation}
l =- \int_{1}^{r}  \frac{\mathrm{d}t }{ \sqrt{1- t^{2n}} }
= \int_{r}^{1}  \frac{\mathrm{d}t }{ \sqrt{1- t^{2n}} }
=\mathrm{arccleaf}_n(r)
\left(= \int_{1}^{r} \mathrm{d}l=[l]^r_1=l(r)-l(1)=l(r) \right)
\label{1.2.5}
\end{equation}

For integrating the left side of the above equation, the initial condition (Eq. (\ref{1.1.2}), $(l,r)=(0,1)$) is applied. The above equation represents the inverse function of the leaf function: cleaf$_n(t)$ \cite{Kaz_cl}. Therefore, the above equation is described as

\begin{equation}
r=\mathrm{cleaf}_n(l) \label{1.2.6}
\end{equation}

Similarly, the function that satisfies Eqs. (\ref{1.1.4})--(\ref{1.1.6}) is defined as sleaf$_n(\overline{l})$.
 In the domain, $0 \leqq \overline{l} \leqq \frac{1}{2} \pi_n$(See Appendix A for constant $\pi_n$), the gradient $d\overline{r}/d\overline{l}$ becomes positive.

\begin{equation}
{\mathrm{d}\overline{l}} = \frac{\mathrm{d}\overline{r} }{ \sqrt{1- \overline{r}^{2n}} } \label{1.2.7}
\end{equation}

The following equation is obtained by integrating the above equation from 0 to $\overline{r}$.

\begin{equation}
\overline{l} = \int_{0}^{\overline{r}}  \frac{\mathrm{d}t }{ \sqrt{1- t^{2n}} }
=\mathrm{arcsleaf}_n(\overline{r})
\left(= \int_{0}^{\overline{r}} \mathrm{d}\overline{l}=[\overline{l}]^{\overline{r}}_0=\overline{l}(\overline{r})-\overline{l}(0)=\overline{l}(\overline{r}) \right)
\label{1.2.8}
\end{equation}

For integrating the left side of the above equation, the initial condition (Eq. (\ref{1.1.5}), $(l,r)=(0,0)$) is applied. The above equation represents the inverse function of the leaf function: sleaf$_n(\overline{l})$ \cite{Kaz_sl}. Therefore, the above equation is described as

\begin{equation}
\overline{r}=\mathrm{sleaf}_n(\overline{l}) \label{1.2.9}
\end{equation}


\subsection{Literature comparison}
Inverse leaf functions based on the basis $n = 1$ represent inverse trigonometric functions.

\begin{equation}
\mathrm{arcsin}(\overline{r}) = \int_{0}^{\overline{r}}  \frac{\mathrm{d}t }{ \sqrt{1- t^2} }
=\mathrm{arcsleaf}_1(\overline{r})=\overline{l}
\qquad (\quad \overline{r}=\mathrm{sin}(\overline{l}) \quad)
\label{1.3.1}
\end{equation}

\begin{equation}
\mathrm{arccos}(r) = \int_{r}^{1}  \frac{\mathrm{d}t }{ \sqrt{1- t^2} }
=\mathrm{arccleaf}_1(r) =l
\qquad (\quad r=\mathrm{cos}(l) \quad)
\label{1.3.2}
\end{equation}

In 1796, Carl Friedrich Gauss presented the lemniscate function\cite{Gauss}. The inverse leaf functions based on the basis $n = 2$ represents inverse functions of the sin and cos lemniscates \cite{roy2017elliptic}.

\begin{equation}
\mathrm{arcsl}(\overline{r}) = \int_{0}^{\overline{r}}  \frac{\mathrm{d}t }{ \sqrt{1- t^4} }
=\mathrm{arcsleaf}_2(\overline{r})=\overline{l} 
\qquad (\quad \overline{r}=\mathrm{sl}(\overline{l}) \quad)
\label{1.3.3}
\end{equation}

\begin{equation}
\mathrm{arccl}(r) = \int_{r}^{1} \frac{\mathrm{d}t }{ \sqrt{1- t^4} }
=\mathrm{arccleaf}_2(r) =l
\qquad (\quad r=\mathrm{cl}(l) \quad)
\label{1.3.4}
\end{equation}

In 1827, Carl Gustav Jacob Jacobi presented the Jacobi elliptic functions \cite{Jacobi}. Compared to Eq. (\ref{1.3.3}), the term $t^2$ is added to the root of the integrand denominator.

\begin{equation}
\mathrm{arcsn}(r,k)=\int_{0}^{r} \frac{\mathrm{d}t}{\sqrt{1-(1+k^2)t^2+k^2 t^4}} \label{1.3.5}
\end{equation}

Eq. (\ref{1.3.5}) represents the inverse Jacobi elliptic function sn, where $k$ is a constant; there are 12 Jacobi elliptic functions, including cn and dn, etc. In Eq. (\ref{1.3.5}), variable $t$ is raised to the fourth power in the denominator. Jacobi did not discuss variable $t$ raised to higher powers as indicated below.

\begin{equation}
\int_{0}^{r} \frac{\mathrm{d}t}{\sqrt{1-t^6}} \: , \: \int_{0}^{r} \frac{\mathrm{d}t}{\sqrt{1-t^8}}\: , \: \int_{0}^{r} \frac{\mathrm{d}t}{\sqrt{1-t^{10} }} \ldots\ \label{1.3.6}
\end{equation}

Thus, historically, the inverse functions have not been discussed in the case of n = 3 or higher \cite{10.2307/108547}\cite{byrd1971handbook} \cite{akhiezerelements} \cite{dixon1894elementary} \cite{lawden1989elliptic}
\cite{mckean1999elliptic}
\cite{walker1996elliptic}
\cite{hancock1958lectures}
\cite{mickens2019generalized}.

\subsection{Originality and Purpose}
A lemniscate is a curve defined by two foci $F_1$ and $F_2$. If the distance between the focal points of $F_1 - F_2$ is 2a (a: constant), then any point $P$ on the lemniscate curve satisfies the equation $PF_1 \cdot PF_2 = a^2$. Jacob Bernoulli first described the lemniscate in 1694 \cite{Bos} \cite{Ayo}. Based on the lemniscate curve, its arc length can be bisected and trisected using a classical ruler and compass \cite{Osler}. Based on this lemniscate, a lemniscate function was proposed by Gauss around the year 1800 \cite{Gauss} \cite{doi:10.1080/00029890.1981.11995279}.
Nishimura proposed a relationship between the product formula for the lemniscate function and Carson's algorithm; it is known as the variant of the arithmetic--geometric mean of Gauss \cite{Nishimura} \cite{CAR}. The Wilker and Huygens-type inequalities have been obtained for Gauss lemniscate functions \cite{doi:10.1080/10652469.2012.684054}.
Deng and Chen established some Shafer--Fink type inequalities for the Gauss
lemniscate function \cite{Deng}. The geometrical characteristics of the lemniscate have been described \cite{Smadja} \cite{Schappacher}.
Mendiratta et al. investigated the geometric properties of functions \cite{Mendiratta}. Levin developed analogs of sine and cosine for the curve to prove the formula \cite{doi:10.1080/00029890.2006.11920331}.
Langer and Singer presented the lemniscate octahedral groups of projective symmetries \cite{Langer}.
As a kinematic control problem, a five body choreography on an algebraic lemniscate was shown as the potential problem for two values of elliptic moduli \cite{LOPEZVIEYRA20191711}. The trajectory generation algorithm was applied by using the shape of the Bernoulli lemmiscate \cite{inproceedings}

Leaf functions are extended lemniscate functions. Various formulas for leaf functions such as the addition theorem of the leaf functions and its application to nonlinear equations have been presented \cite{cmes.2018.02179} \cite{cmes.2019.04472} \cite{cmes.2020.08656}.

In this paper, the geometrical properties of leaf functions for $n = 2$, and the geometric relationship between the angle $\theta$ and lemniscate arc length $l$ are presented using the lemniscate curve. The relations between leaf functions $\mathrm{sleaf}_2(l)$ and $\mathrm{cleaf}_2(l)$ are derived using the geometrical properties of the lemniscate curve, similarity of triangles, and the Pythagorean theorem. In the literature, the relationship equation of $\mathrm{sleaf}_2(l)$ and $\mathrm{cleaf}_2(l)$ is analytically derived; however, it is yet to be derived geometrically \cite{mark}. The relation between $\mathrm{sleaf}_2(l)$ and $\mathrm{cleaf}_2(l)$ can be expressed as

\begin{equation}
(\mathrm{sleaf}_2(l))^2+(\mathrm{cleaf}_2(l))^2
+ (\mathrm{sleaf}_2(l))^2 (\mathrm{cleaf}_2(l))^2 = 1. \label{1.4.1}
\end{equation}

The Eq. (\ref{1.4.1}) was analytically derived. However, it cannot be geometrically derived using the lemniscate curve because
it is not possible to show the geometric relationship of the lemniscate functions $sl (l)$ and $cl (l)$ on a single lemniscate curve. In contrast, phase $l$ of the lemniscate function and angle $\theta$ can be visualized geometrically on a single lemniscate curve. Therefore, in the literature, Eq. (\ref{1.4.1}) is derived using an analytical method without requiring the geometric relationship.

In this paper, the angle $\theta$, phase $l$, and leaf functions $\mathrm{sleaf}_2(l)$ and $\mathrm{cleaf}_2(l)$(or lemniscate functions $sl(l)$ and $cl(l)$) are visualized geometrically on a single lemniscate curve. Eq. (\ref{1.4.1}) is derived based on the geometrical interpretation, similarity of triangles, and Pythagorean theorem.


\section{Geometric relationship with the leaf function $ \mathrm{cleaf}_2(l)$ }

Fig. 1 shows the geometric relationship between the lemniscate curve and $ \mathrm{cleaf}_2(l)$. The $y$ and $x$ axes represent the vertical and horizontal axes, respectively. The equation of the curve is

\begin{equation}
(x^2+y^2)^2=x^2-y^2 \label{2.1}
\end{equation}

If P is an arbitrary point on the lemniscate curve, then the following geometric relation exists.

\begin{equation}
 \mathrm{OP}= \mathrm{cleaf}_2(l) \label{2.2}
\end{equation}

\begin{equation}
 \mathrm{Arc \ \stackrel{\frown}{AP} } = l
 \ (\mathrm{See \ Appendix \ B}).  \label{2.3}
\end{equation}

\begin{equation}
\angle \mathrm{AOP} = \theta \label{2.4}
\end{equation}

When point P is circled along the contour of one leaf, the contour length corresponds to the half cycle $\pi_2$ (See Appendix A for the definition of the constant $\pi_2$). As shown in Fig. 1, with respect to an arbitrary phase $l$, angle $\theta$ must satisfy the following inequality.

\begin{equation}
 \frac{\pi}{4}(4k-1) \leqq \theta \leqq \frac{\pi}{4}(4k+1) \label{2.5}
\end{equation}

Here, $k$ is an integer.

\begin{figure*}[tb]
\begin{center}
\includegraphics[width=0.75 \textwidth]{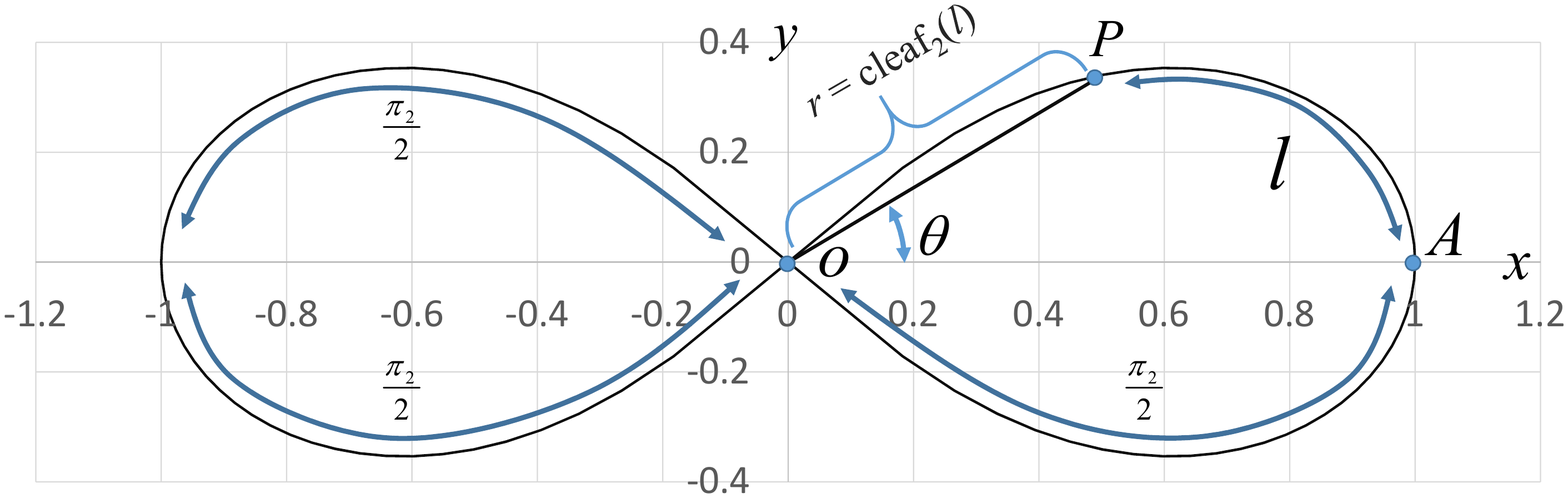}
\caption{Geometric relationship between angle $\theta$ and phase $l$ of the leaf function $\mathrm{cleaf}_n(l)$}
\label{fig1}
\end{center}
\end{figure*}

\section{ Geometric relationship between the trigonometric function and leaf function $\mathrm{cleaf}_2(l)$}

Fig. 2 shows the foci F and F' of the lemniscate curve.
The length of a straight line connecting an arbitrary point $\mathrm{P}$ and one focal point $\mathrm{F}$ is denoted by $\mathrm{PF}$. Similarly, $\mathrm{PF'}$ denotes the length of the line connecting an arbitrary point $\mathrm{P}$ and a second focal point $\mathrm{F'}$. On the curve, the product of $\mathrm{PF}$ and $\mathrm{PF'}$ is constant. The relationship equation is described as \cite{Akopyan}

\begin{equation}
\mathrm{PF} \cdot \mathrm{PF'} = \left( \frac{1}{\sqrt{2}} \right)^2  \label{3.1}
\end{equation}

The coordinates of point P are

\begin{equation}
\mathrm{P}(\mathrm{cleaf}_2(l) \mathrm{cos}(\theta), \mathrm{cleaf}_2(l)\mathrm{sin}(\theta))
=\mathrm{P}(r \mathrm{cos}(\theta), r \mathrm{sin}(\theta))
\label{3.2}
\end{equation}

$\mathrm{PF}$ and $\mathrm{PF'}$ are given by

\begin{equation}
\begin{split}
\mathrm{PF}
&=\sqrt{ \left( \mathrm{cleaf}_2(l) \mathrm{cos}(\theta) - \frac{1}{\sqrt{2}} \right)^2+ \left( \mathrm{cleaf}_2(l) \mathrm{sin}(\theta) \right)^2  } \\
&=\sqrt{ \frac{1}{2} + (\mathrm{cleaf}_2(l))^2
-\sqrt{2} \mathrm{cos}(\theta) \mathrm{cleaf}_2(l)
}
\label{3.3}
\end{split}
\end{equation}

and

\begin{equation}
\begin{split}
\mathrm{PF'}
&=\sqrt{ \left( \mathrm{cleaf}_2(l) \mathrm{cos}(\theta) + \frac{1}{\sqrt{2}} \right)^2+ \left( \mathrm{cleaf}_2(l) \mathrm{sin}(\theta) \right)^2  } \\
&=\sqrt{ \frac{1}{2} + (\mathrm{cleaf}_2(l))^2
+\sqrt{2} \mathrm{cos}(\theta) \mathrm{cleaf}_2(l),
}
\label{3.4}
\end{split}
\end{equation}

respectively.

By substituting Eqs. (\ref{3.3}) and (\ref{3.4}) into Eq. (\ref{3.1}), the relationship equation between the leaf function
$\mathrm{cleaf}_2(l)$ and trigonometric function $\mathrm{cos}(\theta)$ can be derived as

\begin{equation}
\begin{split}
( \mathrm{cleaf}_2(l) )^2=2 (\mathrm{cos}(\theta))^2-1=\mathrm{cos}(2 \theta)
\label{3.5}
\end{split}
\end{equation}

After differentiating Eq.(\ref{3.5}) with respect to $l$,

\begin{equation}
\begin{split}
-2\mathrm{cleaf}_2(l) \sqrt{1-( \mathrm{cleaf}_2(l) )^4}
=-2 \mathrm{sin}(2 \theta) \cdot \frac{\mathrm{d} \theta}{\mathrm{d} l}
\label{3.6}
\end{split}
\end{equation}

The following equation is obtained by combining Eqs. (\ref{3.5}) and (\ref{3.6}).

\begin{equation}
\begin{split}
\frac{ \mathrm{d} \theta}{ \mathrm{d} l}
&=\frac{ \mathrm{cleaf}_2(l) \sqrt{1-( \mathrm{cleaf}_2(l) )^4} }
{ \mathrm{sin}(2 \theta) }
=\frac{ \mathrm{cleaf}_2(l) \sqrt{1-( \mathrm{cleaf}_2(l) )^4} }
{ \sqrt{1-(\mathrm{cos}(2 \theta))^2 } } \\
&=\frac{ \mathrm{cleaf}_2(l) \sqrt{1-( \mathrm{cleaf}_2(l) )^4} }
{ \sqrt{1-( \mathrm{cleaf}_2(l) )^4 } }
=\mathrm{cleaf}_2(l)
\label{3.7}
\end{split}
\end{equation}

The differential equation can be integrated using variable $l$. Parameter $t$ in the integrand is introduced to distinguish it from $l$. The integration of Eq. (\ref{3.7}) in the region $0 \leqq t \leqq l$ yields the following equation (See Appendix C for details).

\begin{equation}
\theta= \int_{0}^{l} \mathrm{cleaf}_2(t) \mathrm{d}t
= \mathrm{arctan} ( \mathrm{sleaf}_2(l) )
\label{3.8}
\end{equation}

Therefore, the following equation holds.

\begin{equation}
\mathrm{tan} (\theta) = \mathrm{sleaf}_2(l)
\label{3.9}
\end{equation}

Eq. (\ref{3.5}) can only be described by variable $l$ as

\begin{equation}
\begin{split}
(\mathrm{cleaf}_2(l))^2
=\mathrm{cos} \left( 2 \int_{0}^{l} \mathrm{cleaf}_2(t) \mathrm{d}t \right)
=\mathrm{cos} (2 \mathrm{arctan} (\mathrm{sleaf}_2(l)))
\label{3.10}
\end{split}
\end{equation}

The phase $\int_{0}^{l} \mathrm{cleaf}_2(t) \mathrm{d}t$ of the cos function is plotted in Fig. 3 through numerical analysis. The horizontal and vertical axes represent variables $l$ and $\theta$, respectively.
As shown in Fig. 3, angle $\theta$ satisfies the
inequality

\begin{equation}
\begin{split}
-\frac{\pi}{4} \leqq
\theta = \int_{0}^{l} \mathrm{cleaf}_2(t) \mathrm{d}t
\leqq \frac{\pi}{4}
\label{3.11}
\end{split}
\end{equation}

Eq. (\ref{3.11}) satisfies the inequality of Eq. (\ref{2.5}) under the condition $k=0$.

Fig. 4 shows the geometric relationship between functions $\mathrm{sleaf_2}(l)$ and $\mathrm{cleaf_2}(l)$. The geometric relation in Eq. (\ref{3.9}) is illustrated in Fig. 4.

Draw a perpendicular line from the point P on the lemniscate to the x-axis. This perpendicular is parallel to the y-axis. Let C be the intersection of this perpendicular and the x-axis. Therefore, the angle $\angle$ OCP is 90 $^{\circ}$.  Next, draw a line perpendicular to the x-axis from the intersection A(1,0) of the lemniscate and the x-axis. Let B be the intersection of this perpendicular and the extension of the straight line OP. Here, $x = \mathrm{OC}$ and $y = \mathrm{CP}$. Substituting these into Eq. (\ref{2.1}) gives

\begin{equation}
(\mathrm{OC}^2+\mathrm{CP}^2)^2
=\mathrm{OC}^2-\mathrm{CP}^2,
\label{3.12}
\end{equation}

\begin{equation}
\angle \mathrm{OCP}=90 ^\circ,
\label{3.12.1}
\end{equation}

and

\begin{equation}
\angle \mathrm{OAB}=90 ^\circ.
\label{3.12.2}
\end{equation}

P and B are moving points, and point A is fixed. When angle $\theta$ is zero, both P and B are at A. The geometric relationship is then expressed as $\mathrm{cleaf}_2(l)=1=$OA and sleaf$_2(l)=0=$AB. As $\theta$ increases, P moves away from A, and it moves along the lemniscate curve. Here, phase $l$ of cleaf$_2(l)$ and sleaf$_2(l)$ corresponds to the length of the arc $\stackrel{\frown}{AP}$. The length of the straight line OP is equal to the value of cleaf$_2(l)$. Point B is the intersection point of the straight lines OP and $x = 1$. In other words, P is the intersection point of the straight line OB and lemniscate curve. As $\theta$ increases, B moves away from A and onto the straight line $x = 1$. That is, it moves in the direction perpendicular to the $x$ axis. The length of straight line AB is equal to the value of sleaf$_2(l)$. When $\theta$ reaches $45^\circ$, P moves to origin O and AB = 1. The length of arc $\stackrel{\frown}{AP}$ (or phase $l$) is $\pi_2/2$. Moreover, cleaf$_2(l)=0=$OP and sleaf$_2(l)=1=$AB.

The relationship $\mathrm{OC:OA=CP:AB}$ is derived by the similarity of triangles  $\mathrm{ \triangle{OAB}\sim\triangle{OCP} }$, as shown in Fig. 4. Thus, the following equation holds.

\begin{equation}
\begin{split}
\mathrm{OC}
&=\frac{ \mathrm{OA} \cdot \mathrm{CP} }{ \mathrm{AB} }
=\frac{ \mathrm{cleaf}_2(l)  \mathrm{sin}(\theta)  }{  \mathrm{sleaf}_2(l)  }
=\frac{ \mathrm{cleaf}_2(l)  \sqrt{1-(\mathrm{cos}(\theta))^2 }  }{  \mathrm{sleaf}_2(l)  } \\
&=\frac{ \mathrm{cleaf}_2(l)  \sqrt{1-\frac{1+(\mathrm{cleaf}_2(l))^2}{2} }  }{  \mathrm{sleaf}_2(l)  }
=\frac{ \mathrm{cleaf}_2(l)  \sqrt{1-(\mathrm{cleaf}_2(l))^2 }  }
{\sqrt{2}  \mathrm{sleaf}_2(l)  }
\label{3.13}
\end{split}
\end{equation}

Eq. (\ref{3.5}) is applied in the transformation process. Similarly, the relationship $\mathrm{OP:PC=OB:BA}$ is derived by the similarity of triangles  $\mathrm{ \triangle{OAB}\sim\triangle{OCP} }$, as shown in Fig . 4. Therefore, the following equation holds.

\begin{equation}
\begin{split}
\mathrm{PC}
&=\frac{ \mathrm{OP} \cdot \mathrm{BA} }{ \mathrm{OB} }
=\frac{ \mathrm{cleaf}_2(l)  \mathrm{sleaf}_2(l)  }{ \sqrt{1+(\mathrm{sleaf}_2(l))^2}  }
\label{3.14}
\end{split}
\end{equation}

By substituting Eqs. (\ref{3.13}) and (\ref{3.14}) into Eq. (\ref{3.12}), the following equation is obtained.

\begin{equation}
\begin{split}
&\left\{ \left(  \frac{ \mathrm{cleaf}_2(l)  \sqrt{1-(\mathrm{cleaf}_2(l))^2 }  }
{\sqrt{2}  \mathrm{sleaf}_2(l)  } \right)^2 +
\left( \frac{ \mathrm{cleaf}_2(l)  \mathrm{sleaf}_2(l)  }
{ \sqrt{1+(\mathrm{sleaf}_2(l))^2}  } \right)^2  \right\}^2 \\
&=\left(  \frac{ \mathrm{cleaf}_2(l)  \sqrt{1-(\mathrm{cleaf}_2(l))^2 }  }
{\sqrt{2}  \mathrm{sleaf}_2(l)  } \right)^2 -
\left( \frac{ \mathrm{cleaf}_2(l)  \mathrm{sleaf}_2(l)  }
{ \sqrt{1+(\mathrm{sleaf}_2(l))^2}  } \right)^2
\label{3.15}
\end{split}
\end{equation}

By rearranging Eq. (\ref{3.15}), the following equation is obtained.

\begin{equation}
\begin{split}
\frac{ (\mathrm{cleaf}_2(l))^2  \{ -1+ (\mathrm{sleaf}_2(l))^2
+(\mathrm{cleaf}_2(l))^2+(\mathrm{sleaf}_2(l))^2(\mathrm{cleaf}_2(l))^2 \}
 \{ \cdots \} } {4 (\mathrm{sleaf}_2(l))^2 \{1+ (\mathrm{sleaf}_2(l))^2 \}^2 }=0
\label{3.16}
\end{split}
\end{equation}

\begin{equation}
\begin{split}
& \{ \cdots \} =2(\mathrm{sleaf}_2(l))^2+6 (\mathrm{sleaf}_2(l))^4
+4(\mathrm{sleaf}_2(l))^6 \\
&- (\mathrm{cleaf}_2(l))^2
\{ 1+ 3\mathrm{sleaf}_2(l))^2+4\mathrm{sleaf}_2(l))^4 \}
+ (\mathrm{cleaf}_2(l))^4 \{ 1+ \mathrm{sleaf}_2(l))^2 \}
\label{3.16}
\end{split}
\end{equation}

For arbitrary $l$, $\mathrm{cleaf}_2(l) \neq 0$ and $\{ \cdots \} \neq 0$. The relationship between $\mathrm{sleaf}_2(l)$ and $\mathrm{cleaf}_2(l)$ can then be obtained as Eq. (\ref{1.4.1}).

\begin{figure*}[tb]
\begin{center}
\includegraphics[width=0.75 \textwidth]{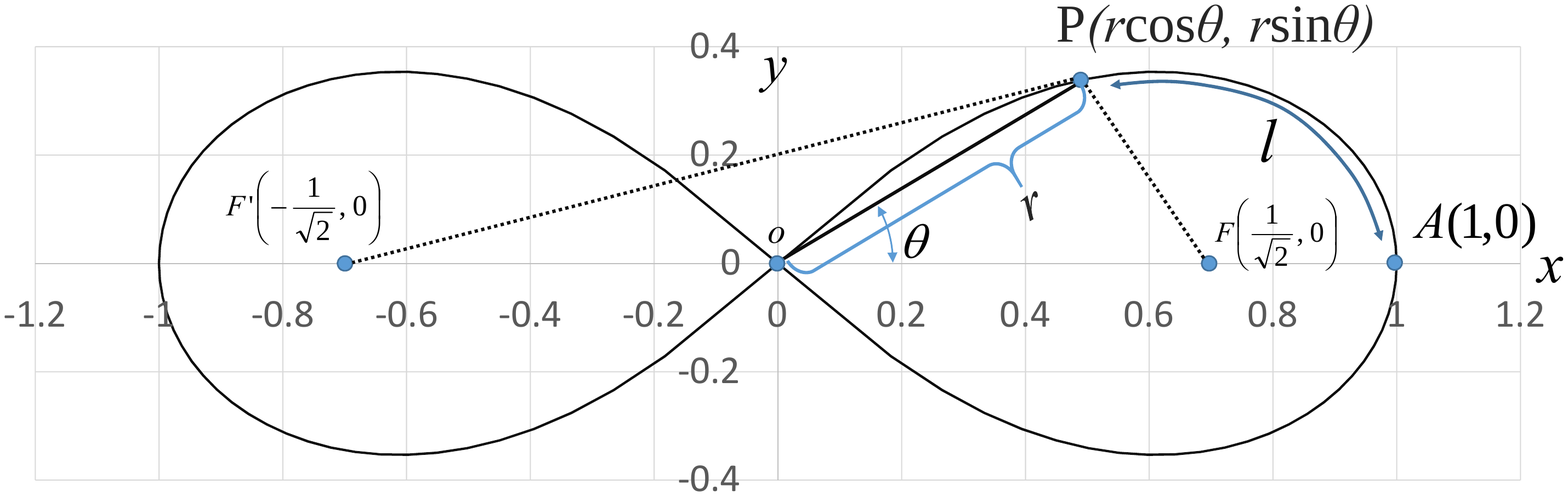}
\caption{Lemniscate focus}
\label{fig2}
\end{center}
\end{figure*}

\begin{figure*}[tb]
\begin{center}
\includegraphics[width=0.75 \textwidth]{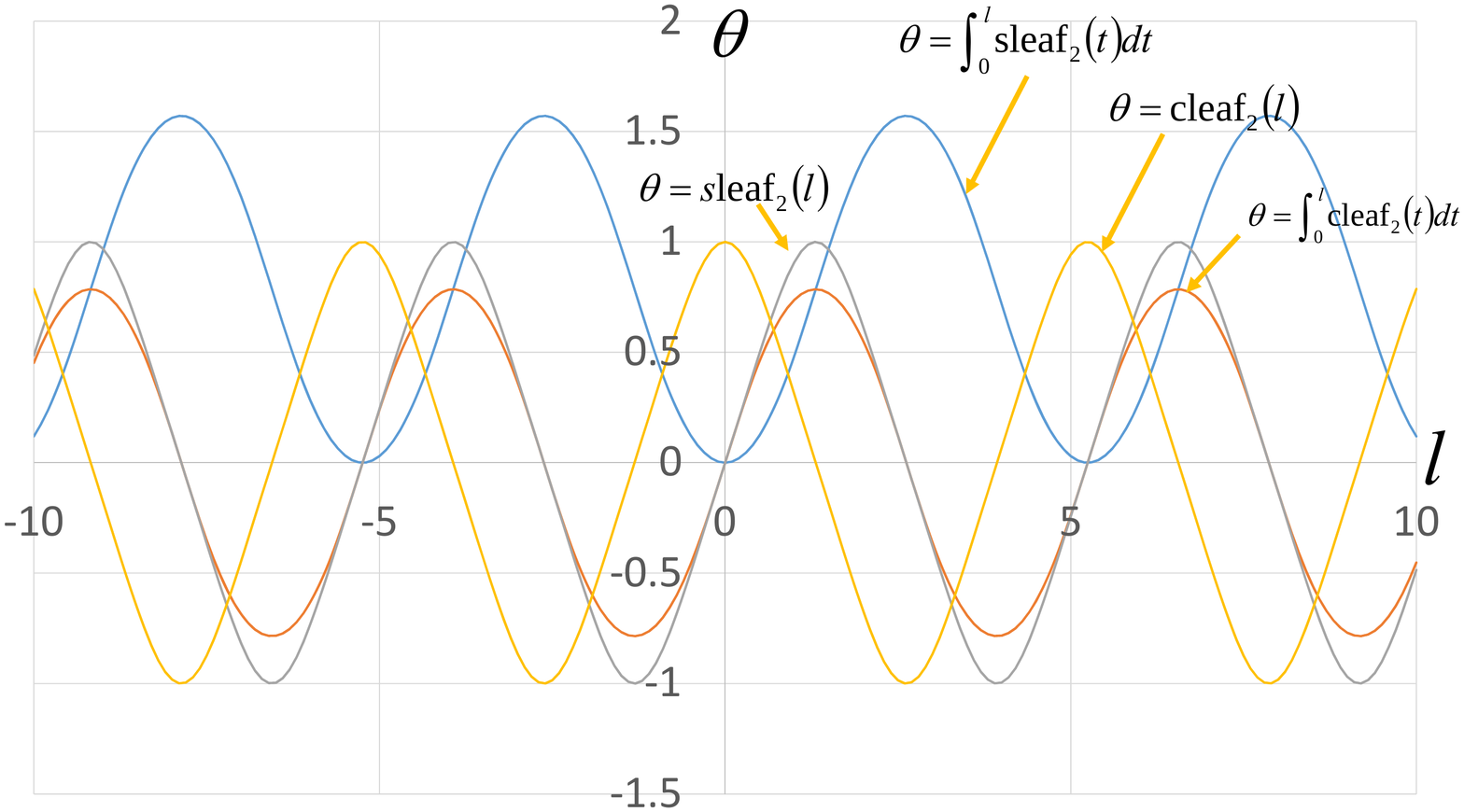}
\caption{ Curves of leaf functions ($\mathrm{sleaf}_2(l)$  and $\mathrm{cleaf}_2(l)$) and the integrated leaf functions ($\int_{0}^{l} \mathrm{sleaf}_2(t) \mathrm{d}t$  and $\int_{0}^{l} \mathrm{cleaf}_2(t) \mathrm{d}t$ )}
\label{fig3}
\end{center}
\end{figure*}

\begin{figure*}[tb]
\begin{center}
\includegraphics[width=0.75 \textwidth]{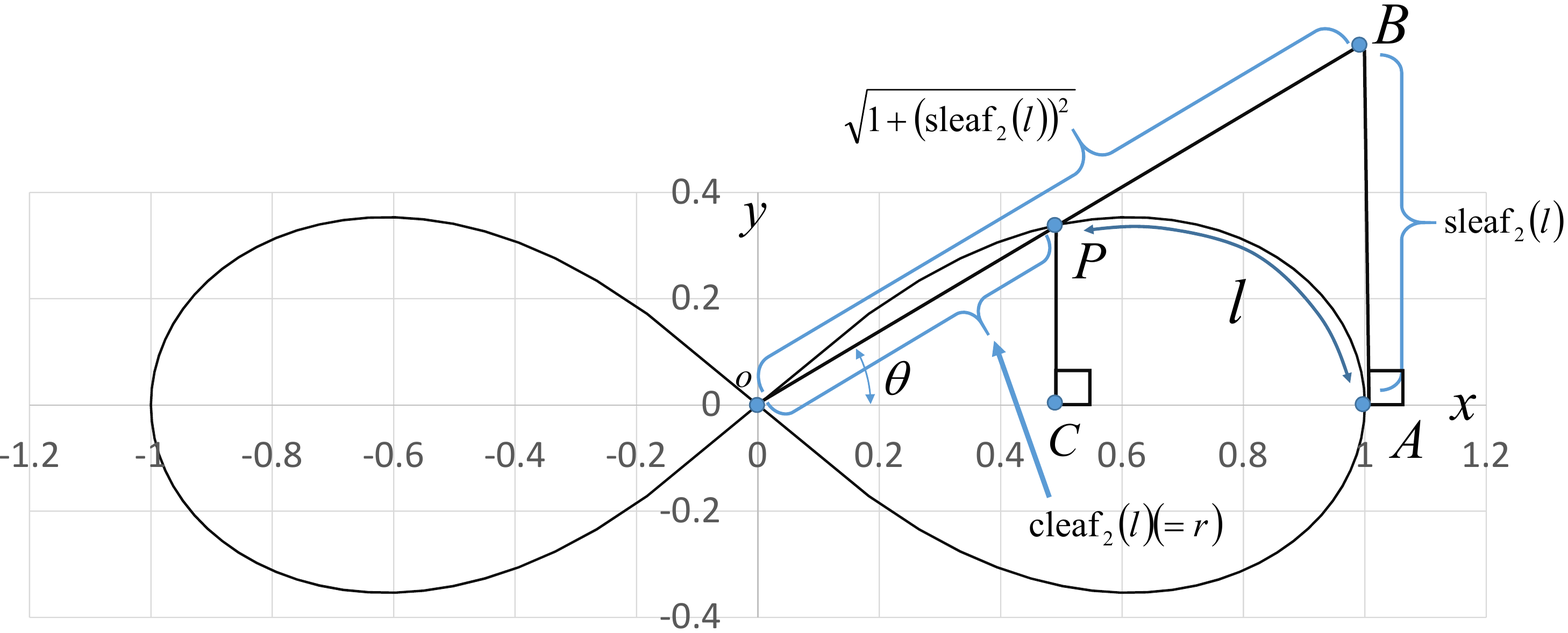}
\caption{Geometric relationship between leaf functions $\mathrm{sleaf}_2(l)$  and  $\mathrm{cleaf}_2(l)$ }
\label{fig4}
\end{center}
\end{figure*}

\section{ Geometric relationship of the leaf function $\mathrm{sleaf}_2( \overline{l})$ }

Fig. 5 shows the geometric relationship between length $\mathrm{sleaf}_2(\overline{l})$ and lemniscate curve inclined at $45^\circ$. In Fig. 5, the $y$ and $x$ axes represent the vertical and horizontal axes, respectively.  The equation of this curve is given as

\begin{equation}
(x^2+y^2)^2=2 x y  \label{4.1}
\end{equation}

If $\overline{\mathrm{P}}$ is an arbitrary point on the lemniscate curve, the following geometric relation exists.

\begin{equation}
  \mathrm{O\overline{P} }= \mathrm{sleaf}_2( \overline{l})  \label{4.2}
\end{equation}

\begin{equation}
 \mathrm{Arc \ O\overline{P} } = \overline{l}
\ (\mathrm{See \ Appendix D} ) \label{4.3}
\end{equation}

\begin{equation}
\angle \mathrm{\overline{D}O\overline{P}} = \overline{\theta}  \label{4.4}
\end{equation}

In Fig. 5, for an arbitrary variable $\overline{l}$, the range of angle $\overline{\theta}$ is given by

\begin{equation}
 k \pi \leqq \overline{\theta} \leqq \frac{\pi}{2}(2k+1).   \label{4.5}
\end{equation}

Here, $k$ is an integer.

\begin{figure*}[tb]
\begin{center}
\includegraphics[width=0.75 \textwidth]{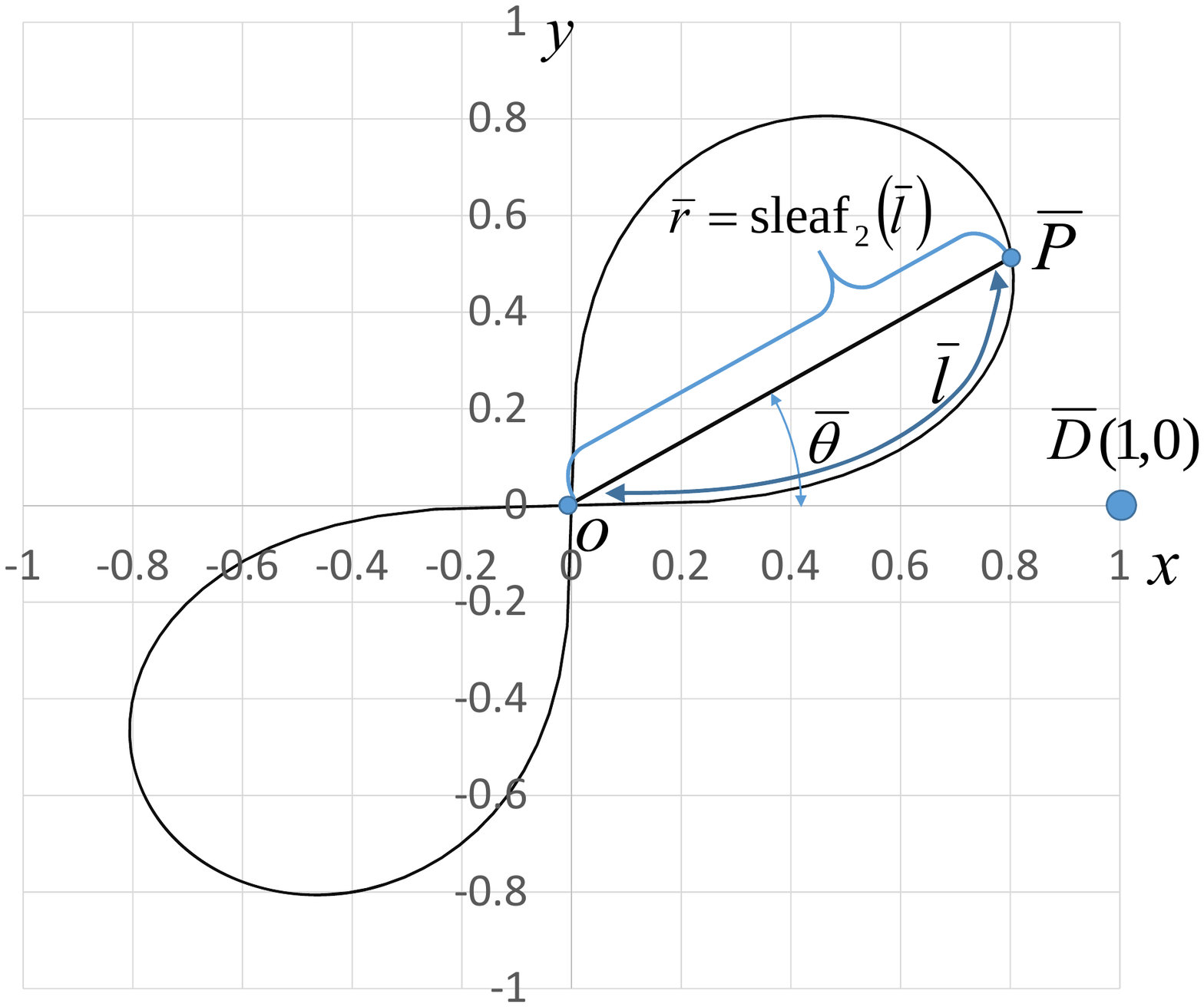}
\caption{ Geometric relationship between angle $\overline{\theta}$ and phase $\overline{l}$ of leaf function $\mathrm{sleaf}_2(\overline{l})$}
\label{fig5}
\end{center}
\end{figure*}

\section{ Geometric relationship between trigonometric function and leaf function $\mathrm{sleaf}_2(l)$}

Fig. 6 shows foci $\overline{\mathrm{F}}$ and $\overline{\mathrm{F'}}$ of the lemniscate curve inclined at an angle of $45^\circ$. This curve has the same relation equation as shown in Fig. 2.

\begin{equation}
\mathrm{\overline{PF} \cdot \overline{PF'}} = \left( \frac{1}{\sqrt{2}} \right)^2   \label{5.1}
\end{equation}

\begin{figure*}[tb]
\begin{center}
\includegraphics[width=0.75 \textwidth]{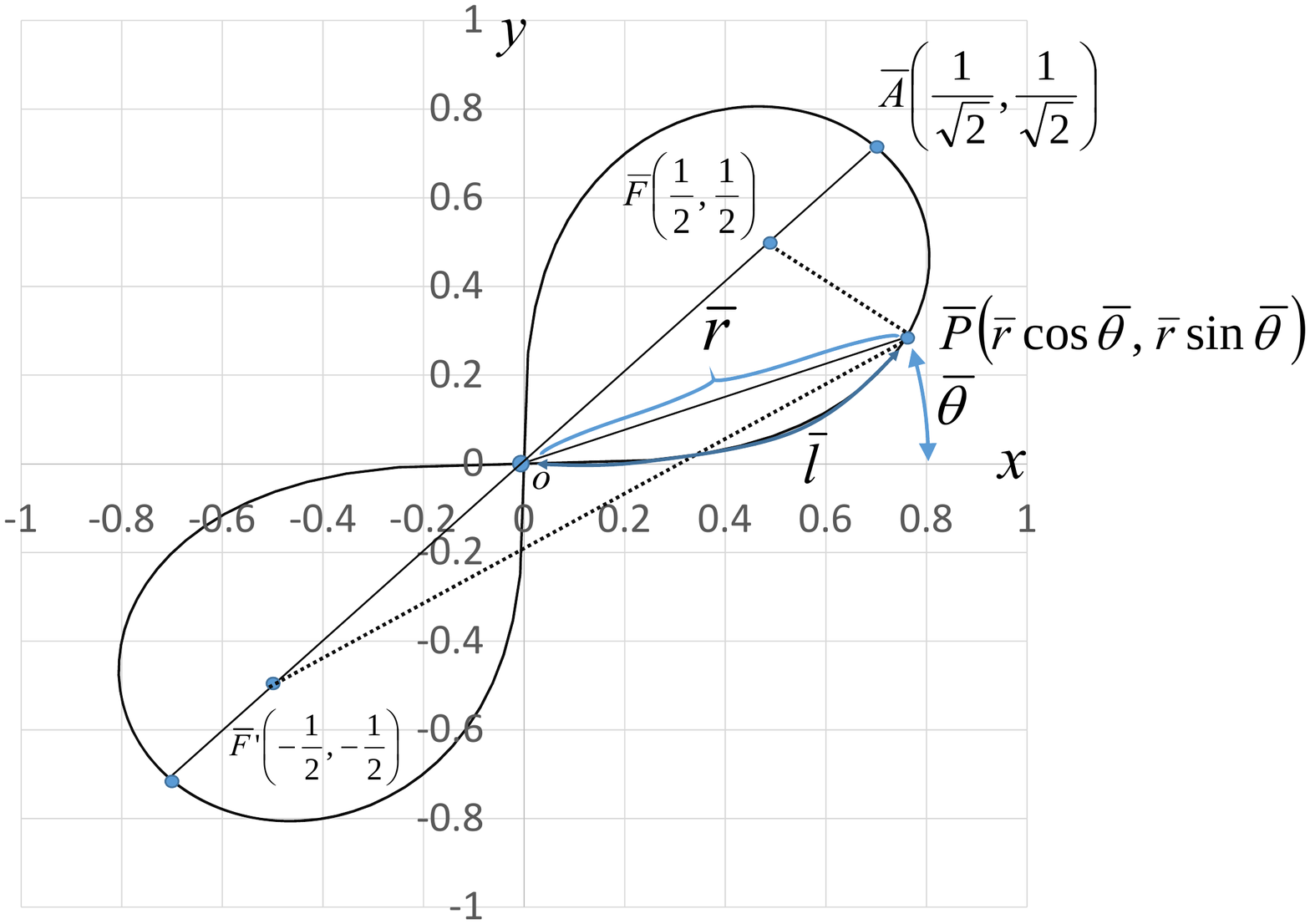}
\caption{ Lemniscate curve inclined at  $45^\circ$ }
\label{fig6}
\end{center}
\end{figure*}

The coordinates of point $\overline{\mathrm{P}}$ on the lemniscate curve inclined at an angle of $45^\circ$ are given as

\begin{equation}
\overline{\mathrm{P}} (\mathrm{sleaf}_2(\overline{l}) \mathrm{cos}(\overline{\theta}), \mathrm{sleaf}_2(\overline{l})\mathrm{sin}(\overline{\theta}))
=\overline{\mathrm{P}} (\overline{r} \mathrm{cos}(\overline{\theta}), \overline{r} \mathrm{sin}(\overline{\theta}))
\label{5.2}
\end{equation}

Lengths $\overline{\mathrm{PF}}$ and $\overline{\mathrm{PF'}}$ are expressed by

\begin{equation}
\begin{split}
\overline{\mathrm{PF}}
&=\sqrt{ \left( \mathrm{sleaf}_2(\overline{l}) \mathrm{cos}(\overline{\theta}) - \frac{1}{ 2 } \right)^2
+ \left( \mathrm{sleaf}_2(\overline{l}) \mathrm{sin}(\overline{\theta})  - \frac{1}{ 2 } \right)^2   } \\
&=\sqrt{ \frac{1}{2} + (\mathrm{sleaf}_2(\overline{l}))^2
-( \mathrm{sin}(\overline{\theta}) +  \mathrm{cos}(\overline{\theta})) \mathrm{sleaf}_2(\overline{l})
}
\label{5.3}
\end{split}
\end{equation}

and

\begin{equation}
\begin{split}
\overline{\mathrm{PF'}}
&=\sqrt{ \left( \mathrm{sleaf}_2(\overline{l}) \mathrm{cos}(\overline{\theta})
+ \frac{1}{ 2 } \right)^2
+ \left( \mathrm{sleaf}_2(\overline{l}) \mathrm{cos}(\overline{\theta})  + \frac{1}{ 2 } \right)^2   } \\
&=\sqrt{ \frac{1}{2} + (\mathrm{sleaf}_2(\overline{l}))^2
+( \mathrm{sin}(\overline{\theta}) +  \mathrm{cos}(\overline{\theta})) \mathrm{sleaf}_2(\overline{l}).
}
\label{5.4}
\end{split}
\end{equation}

By substituting Eqs. (\ref{5.3}) and (\ref{5.4}) into Eq. (\ref{5.1}), the relationship equation between the leaf function
$\mathrm{sleaf}_2(\overline{l})$ and the trigonometric function $\mathrm{sin}( \overline{\theta})$ can be derived as

\begin{equation}
\begin{split}
( \mathrm{sleaf}_2(\overline{l}) )^2=2 \mathrm{sin}(\overline{\theta}) \mathrm{cos}(\overline{\theta})  = \mathrm{sin}(2\overline{\theta})
\label{5.5}
\end{split}
\end{equation}

The following equation is obtained by differentiating Eq. (\ref{5.5}) with respect to the variable $\overline{l}$.

\begin{equation}
\begin{split}
2\mathrm{sleaf}_2(\overline{l}) \sqrt{1-( \mathrm{sleaf}_2(\overline{l}) )^4}
=2 \mathrm{cos}(2 \overline{\theta}) \cdot \frac{\mathrm{d} \overline{\theta}}{\mathrm{d} \overline{l}}
\label{5.6}
\end{split}
\end{equation}

After applying Eq. (\ref{5.5}), the equation is transformed as

\begin{equation}
\begin{split}
\frac{ \mathrm{d} \overline{\theta} }{ \mathrm{d} \overline{l} }
&=\frac{ \mathrm{sleaf}_2( \overline{l} )  \sqrt{1-( \mathrm{sleaf}_2( \overline{l}  ) )^4} }
{ \mathrm{cos}(2 \overline{\theta} )}
=\frac{ \mathrm{sleaf}_2( \overline{l} ) \sqrt{1-( \mathrm{sleaf}_2(  \overline{l}  ) )^4} }
{ \sqrt{1-(\mathrm{sin}(2 \overline{\theta}  ))^2 } } \\
&=\frac{ \mathrm{sleaf}_2(\overline{l} ) \sqrt{1-( \mathrm{sleaf}_2( \overline{l} ) )^4} }
{ \sqrt{1-( \mathrm{sleaf}_2(\overline{l} ) )^4 } }
=\mathrm{sleaf}_2(\overline{l} )
\label{5.7}
\end{split}
\end{equation}

The differential equation is integrated by variable $\overline{l}$. Parameter $t$ is introduced to distinguish the parameter from variable $\overline{l}$ in the integration region. The integration of Eq. (\ref{5.7}) in region  $0 \leqq  t \leqq  \overline{l}$ (See Appendix C) yields

\begin{equation}
\overline{\theta}= \int_{0}^{\overline{l} } \mathrm{sleaf}_2(t) \mathrm{d}t
= - \mathrm{arctan} ( \mathrm{cleaf}_2(\overline{l}) ) +\frac{\pi}{4},
\label{5.8}
\end{equation}

and the following equation holds.

\begin{equation}
\mathrm{tan} \left( \frac{\pi}{4} - \overline{\theta} \right) = \mathrm{cleaf}_2(\overline{l}).
\label{5.9}
\end{equation}

Using Eq. (\ref{5.8}), Eq. (\ref{5.5}) can be described by variable $\overline{l}$ as

\begin{equation}
\begin{split}
(\mathrm{sleaf}_2( \overline{l} ))^2
=\mathrm{sin} \left( 2 \int_{0}^{\overline{l}} \mathrm{sleaf}_2(t) \mathrm{d}t  \right)
=\mathrm{sin} \left(2 ( \frac{\pi}{4} - \mathrm{arctan} (\mathrm{cleaf}_2(\overline{l}) )  ) \right)
=\mathrm{cos} \left(2 \mathrm{arctan} (\mathrm{cleaf}_2(\overline{l}) )  ) \right)
\label{5.10}
\end{split}
\end{equation}

The curve of phase $ \overline{\theta} = \int_{0}^{ \overline{l} } \mathrm{sleaf}_2(t) \mathrm{d}t$  is plotted in Fig. 3 through numerical analysis.  The horizontal and vertical axes represent variables $\overline{l}$ and $\overline{\theta}$, respectively.  As shown in Fig. 3, angle $\overline{\theta}$ satisfies the following range.

\begin{equation}
\begin{split}
0 \leqq \overline{\theta} = \int_{0}^{\overline{l}} \mathrm{sleaf}_2(t) \mathrm{d}t
\leqq \frac{\pi}{2}
\label{5.11}
\end{split}
\end{equation}

\begin{figure*}[tb]
\begin{center}
\includegraphics[width=0.75 \textwidth]{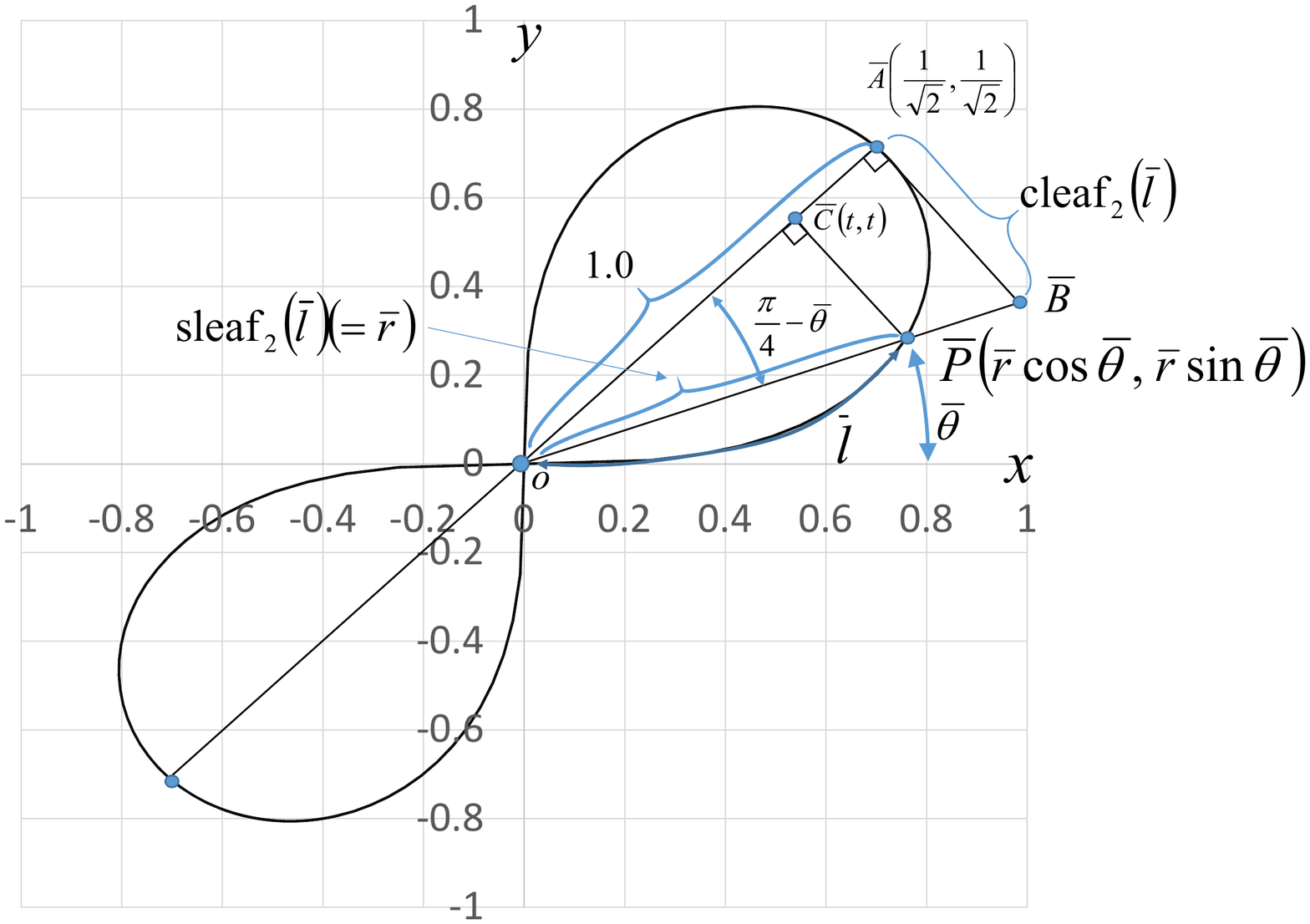}
\caption{ Geometric relationship based on the lemniscate curve inclined at $45^\circ$}
\label{fig7}
\end{center}
\end{figure*}

Eq. (\ref{5.11}) satisfies the range of Eq. (\ref{4.5}) under the condition $k=0$ . Fig. 7 shows the lemniscate curve inclined at $45^\circ$. The geometric relation of Eq. (\ref{5.9}) is added to Fig. 7.

\begin{equation}
\begin{split}
\angle \mathrm{O\overline{CP}}=90 ^\circ
\label{5.11_1}
\end{split}
\end{equation}
\begin{equation}
\begin{split}
\angle \mathrm{O\overline{AB}}=90 ^\circ
\label{5.11_2}
\end{split}
\end{equation}

Let $\overline{C}(t,t)$ be the coordinates on line $\mathrm{O}\overline{\mathrm{A}}$, as shown in Fig. 7. Points $\overline{\mathrm{P}}$ and $\overline{\mathrm{B}}$ are moving points, and point $\overline{\mathrm{A}}$ is fixed. When angle $\overline{\theta}$ is zero, $\overline{\mathrm{P}}$ is at origin O, and $\overline{\mathrm{B}}$ is on the $x$-axis at $(x,y)=(\sqrt{2}$,0). The geometric relationship is described as $\mathrm{cleaf}_2(\overline{l})=1=\overline{\mathrm{AB}}$ and sleaf$_2(\overline{l})=0$. As $\overline{\theta}$ increases, $\overline{\mathrm{P}}$ moves away from origin O and along the lemniscate curve. The phase $\overline{l}$ of both cleaf$_2(\overline{l})$ and sleaf$_2(\overline{l})$ corresponds to the length of arc $\stackrel{\frown}{\mathrm{OP}}$.  The length of the straight line $\overline{\mathrm{OP}}$ is equal to the value of sleaf$_2(\overline{l}) (=\overline{r})$.  Point $\overline{\mathrm{B}}$ is the intersection point of straight lines $\overline{\mathrm{OP}}$ and $y=-x+\sqrt{2}$. In other words, $\overline{\mathrm{P}}$ is the intersection point of the straight line $\overline{\mathrm{OB}}$ and the lemniscate curve. As $\overline{\theta}$ increases, $\overline{\mathrm{B}}$ moves away from point $(x,y)=(\sqrt{2}$,0) on a straight line $y=-x+\sqrt{2}$.  Here, the length of the straight line $\overline{\mathrm{AB}}$ is equal to the value of cleaf$_2(\overline{l})$. When $\overline{\theta}$ reaches $45^\circ$, $\overline{\mathrm{P}}$ moves to point $\overline{\mathrm{A}}$. The length of arc $\stackrel{\frown}{\mathrm{OA}}$ (or phase $\overline{l}$) becomes $\pi_2/2$. Furthermore, cleaf$_2(\overline{l})=0$ and sleaf$_2(\overline{l})=1=\mathrm{O\overline{A}}$. The linear equation $\overline{\mathrm{CP}}$ is given by

\begin{equation}
y=-x+2t.
\label{5.12}
\end{equation}

By substituting Eq. (\ref{5.12}) into Eq. (\ref{4.1}) and solving for variable $x$, four solutions can be obtained as

\begin{equation}
x= \frac{1}{2} \left\{ 2t - \sqrt{ -1-4t^2-\sqrt{1+16t^2} }    \right\}
\label{5.13}
\end{equation}

\begin{equation}
x= \frac{1}{2} \left\{ 2t + \sqrt{ -1-4t^2-\sqrt{1+16t^2} }    \right\}
\label{5.14}
\end{equation}

\begin{equation}
x= \frac{1}{2} \left\{ 2t - \sqrt{ -1-4t^2+\sqrt{1+16t^2} }    \right\}
\label{5.15}
\end{equation}

\begin{equation}
x= \frac{1}{2} \left\{ 2t + \sqrt{ -1-4t^2+\sqrt{1+16t^2} }    \right\}
\label{5.16}
\end{equation}

As Eqs. (\ref{5.13}) and (\ref{5.14}) include imaginary numbers, the solutions for $x$ using both Eqs. (\ref{5.15}) and (\ref{5.16}) are determined by the intersection points of line $\overline{\mathrm{CP}}$ and the lemniscate curve, as shown in Fig. 7. The larger $x$ value is given by Eq. (\ref{5.16}).  That is, the coordinates of point $\overline{\mathrm{P}}$ can be expressed as

\begin{equation}
\overline{\mathrm{P}} \left( t+  \frac{1}{2}  \sqrt{ -1-4t^2+\sqrt{1+16t^2} }, \quad
 t-  \frac{1}{2}  \sqrt{ -1-4t^2+\sqrt{1+16t^2} }     \right)
\label{5.17}
\end{equation}

Therefore, the length $\overline{\mathrm{CP}}$ is expressed as

\begin{equation}
\overline{\mathrm{CP}} = \frac{1}{\sqrt{2}}  \sqrt{ -1-4t^2+\sqrt{1+16t^2} }.
\label{5.18}
\end{equation}

The following equation is obtained from the Pythagorean theorem of the triangle
$\triangle{ \mathrm{O\overline{PC}} }$.

\begin{equation}
\mathrm{O\overline{P}^2 =\overline{CP}^2 +O\overline{C}^2  }
\label{5.19}
\end{equation}

Substitution of Eqs. (\ref{4.2}) and (\ref{5.18}) into Eq. (\ref{5.19}) yields

\begin{equation}
(\mathrm{sleaf}_2(t))^2 = \frac{1}{2}  ( -1-4t^2+\sqrt{1+16t^2} ) + 2t^2.
\label{5.20}
\end{equation}

The length ratio is $\mathrm{ O\overline{C}:O\overline{A}=\overline{CP}:\overline{AB} }$ owing to the similarity of triangles $\mathrm{ \triangle{O\overline{AB}}\sim\triangle{O\overline{CP}} }$.

\begin{equation}
\sqrt{2}t:1= \frac{1}{\sqrt{2}}  \sqrt{ -1-4t^2+\sqrt{1+16t^2} } : \mathrm{cleaf}_2(t)
\label{5.21}
\end{equation}

The elimination of variable $t$ from Eqs. (\ref{5.20}) and (\ref{5.21}) yields the relation equation, Eq. (\ref{1.4.1}).


\section{Numerical results}
\begin{table}
\begin{center}
\caption{ Numerical data of arc length $l$, angle $\theta$, and leaf functions sleaf$_2(l)$ and cleaf$_2(l)$ for Fig. 4. }
\label{tab1}
\begin{tabular}{ccccc}
\hline\noalign{\smallskip}
$l$ & $\theta$(radian)  $ \left(  =\int_{0}^{l} \mathrm{cleaf}_2(t) \mathrm{d}t  \right) $   & $\theta$(degree) & sleaf$_2(l)$ (=tan $\theta$) & cleaf$_2(l) (=r)$  \\
\noalign{\smallskip}\hline\noalign{\smallskip}
0.0&	0.00000$\cdots$&	0.00000$\cdots$&	0.00000$\cdots$&	1.00000$\cdots$\\
0.1&	0.09966$\cdots$&	 5.7105$\cdots$&	0.09999$\cdots$&	0.99004$\cdots$\\
0.2&	0.19736$\cdots$&	11.3081$\cdots$&	0.19996$\cdots$&	0.96078$\cdots$\\
0.3&	0.29123$\cdots$&	16.6864$\cdots$&	0.29975$\cdots$&	0.91384$\cdots$\\
0.4&	0.37962$\cdots$&	21.7509$\cdots$&	0.39897$\cdots$&	0.85167$\cdots$\\
0.5&	0.46115$\cdots$&	26.4223$\cdots$&	0.49689$\cdots$&	0.77715$\cdots$\\
0.6&	0.53474$\cdots$&	30.6385$\cdots$&	0.59230$\cdots$&	0.69323$\cdots$\\
0.7&	0.59958$\cdots$&	34.3534$\cdots$&	0.68352$\cdots$&	0.60260$\cdots$\\
0.8&	0.65511$\cdots$&	37.5355$\cdots$&	0.76831$\cdots$&	0.50756$\cdots$\\
0.9&	0.70100$\cdots$&	40.1646$\cdots$&	0.84400$\cdots$&	0.40985$\cdots$\\
1.0&	0.73704$\cdots$&	42.2294$\cdots$&	0.90768$\cdots$&	0.31073$\cdots$\\
1.1&	0.76313$\cdots$&	43.7243$\cdots$&	0.95643$\cdots$&	0.21098$\cdots$\\
1.2&	0.77923$\cdots$&	44.6468$\cdots$&	0.98774$\cdots$&	0.11102$\cdots$\\
1.3&	0.78533$\cdots$&	44.9965$\cdots$&	0.99987$\cdots$&	0.01102$\cdots$\\
\noalign{\smallskip}\hline
\end{tabular}
\end{center}
\end{table}
Table 1 lists the numerical values for Fig. 4.
The angle $\theta$ and values of the leaf function (sleaf$_2(l)$ and cleaf$_2(l)$ ) are calculated along the arc length $l$. The values of leaf functions cleaf$_2(l)$ are calculated by Eqs. (\ref{1.1.1})--(\ref{1.1.3}). The values of leaf functions sleaf$_2(l)$ are calculated by Eq. (\ref{1.1.4})--(\ref{1.1.6}). Based on these data, the numerical data of functions sleaf$_2(l)$ and cleaf$_2(l)$ can be confirmed using Eq. (\ref{1.4.1}). The function cleaf$_2(l) (=r)$ can also be confirmed by using Eq. (\ref{B7}). Angle $\theta$ can be calculated by using Eq. (\ref{3.8}).

Table 2 shows the numerical values for Fig. 7. The angle $\overline{\theta}$ and the values of leaf function (sleaf$_2(\overline{l})$ and cleaf$_2(\overline{l})$ ) are calculated along the arc length $\overline{l}$. Based on these data, the function sleaf$_2(\overline{l}) (=\overline{r})$ can also be confirmed using Eq. (\ref{D7}). The angle $\overline{\theta}$ can be calculated using Eq. (\ref{5.8}).

\begin{table}
\begin{center}
\caption{ Numerical data of arc length $\overline{l}$, angle $\overline{\theta}$, and leaf functions sleaf$_2(\overline{l})$ and cleaf$_2(\overline{l})$ for Fig. 7. }
\label{tab1}
\begin{tabular}{ccccc}
\hline\noalign{\smallskip}
$\overline{l}$ & $\overline{\theta}$(radian)  $ \left(  =\int_{0}^{\overline{l}} \mathrm{sleaf}_2(t) \mathrm{d}t  \right) $   & $\overline{\theta}$(degree) & sleaf$_2(\overline{l})$ $(=\overline{r})$ & cleaf$_2(\overline{l})$  (=tan $\overline{\theta}$)  \\
\noalign{\smallskip}\hline\noalign{\smallskip}
0.0&	0.00000$\cdots$&0.0000$\cdots$&	0.00000$\cdots$&	1.00000$\cdots$\\
0.1&	0.00499$\cdots$&0.2864$\cdots$&	0.09999$\cdots$&	0.99004$\cdots$\\
0.2&	0.01999$\cdots$&1.1458$\cdots$&	0.19996$\cdots$&	0.96078$\cdots$\\
0.3&	0.04498$\cdots$&2.5776$\cdots$&	0.29975$\cdots$&	0.91384$\cdots$\\
0.4&	0.07993$\cdots$&4.5797$\cdots$&	0.39897$\cdots$&	0.85167$\cdots$\\
0.5&	0.12474$\cdots$&7.1470$\cdots$&	0.49689$\cdots$&	0.77715$\cdots$\\
0.6&	0.17922$\cdots$&10.2689$\cdots$&	0.59230$\cdots$&	0.69323$\cdots$\\
0.7&	0.24306$\cdots$&13.9264$\cdots$&	0.68352$\cdots$&	0.60260$\cdots$\\
0.8&	0.31571$\cdots$&18.0893$\cdots$&	0.76831$\cdots$&	0.50756$\cdots$\\
0.9&	0.39642$\cdots$&22.7133$\cdots$&	0.84400$\cdots$&	0.40985$\cdots$\\
1.0& 0.48411$\cdots$&27.7379$\cdots$&	0.90768$\cdots$&	0.31073$\cdots$\\
1.1&	0.57746$\cdots$&33.0860$\cdots$&	0.95643$\cdots$&	0.21098$\cdots$\\
1.2&	0.67482$\cdots$&38.6645$\cdots$&	0.98774$\cdots$&	0.11102$\cdots$\\
1.3&	0.77436$\cdots$&44.3681$\cdots$&	0.99987$\cdots$&	0.01102$\cdots$\\
\noalign{\smallskip}\hline
\end{tabular}
\end{center}
\end{table}

\section{Conclusion}
\label{Conclusion}
Based on the geometric properties of the lemniscate curve, the geometric relationship among angle $\theta$, lemniscate length $l$, and leaf functions $\mathrm{ sleaf}_2(l)$ and $\mathrm{ cleaf}_2 (l)$ were shown on the lemniscate curve. Using the similarity of triangles and the Pythagorean theorem,  the relationship equation of leaf functions $\mathrm{ sleaf}_2(l) $ and $\mathrm{ cleaf}_2 (l)$ was derived.

\nocite{*}
\bibliographystyle{unsrt}  
\bibliography{references}  

\appendix
\def\thesection{Appendix}
\section{A}
\label{A}

The symbols $\pi_n$ are a constant given by

\begin{equation}
\pi_n=2 \int_{0}^{1} \frac{\mathrm{d}t}
{\sqrt{1-t^{2n}}} \ (n=1,2,3 \cdots).
\label{A1}
\end{equation}

The numerical data of the symbol $\pi_n$ are summarized in the table 3.

\begin{table}
\begin{center}
\caption{ Values of constants $\pi_n$ }
\label{tab3}
\begin{tabular}{ll}
\hline\noalign{\smallskip}
$n$ & $\pi_n$  \\
\noalign{\smallskip}\hline\noalign{\smallskip}
1	&	3.1415 $\cdots$	\\
2	&	2.6220 $\cdots$	\\
3	&	2.4286 $\cdots$	\\
$\cdots$	&	 $\cdots$	\\
\noalign{\smallskip}\hline
\end{tabular}
\end{center}
\end{table}

\appendix
\def\thesection{Appendix}
\section{B}
\label{B}
The equation (\ref{2.1}) of Cartesian coordinate system is transformed to the following equation of   polar coordinates \cite{Kaz_cl} \cite{Cox}.

\begin{equation}
r^2=\mathrm{cos}(2\theta) \label{B1}
\end{equation}

The variable $r$ and $\theta$ represents OP and $\angle$AOP in Fig. 1. The arc length in the cartesian and polar coordinates is given by

\begin{equation}
\sqrt{dx^2+dy^2}=\sqrt{dr^2+r^2 d\theta^2} \label{B2}
\end{equation}

The arc length $l$ of the lemniscate with polar coordinates is given by

\begin{equation}
l= \int_{r}^{1} \sqrt{1+r^2 \left( \frac{ \mathrm{d} \theta}{ \mathrm{d} r} \right)^2 }  \mathrm{d}r
\label{B3}
\end{equation}

By differentiating Eq. (\ref{B1}) with respect to variable $\theta$,

\begin{equation}
2r \frac{\mathrm{d} r}{\mathrm{d} \theta} =-2 \mathrm{sin}(2\theta)
\label{B4}
\end{equation}

Then, by applying Eq (\ref{B4}),

\begin{equation}
1+r^2 \left( \frac{ \mathrm{d} \theta}{ \mathrm{d} r} \right)^2=1+r^2 \left(- \frac{r}{\mathrm{sin}(2\theta)} \right)^2 =1+\frac{r^4}{1-r^4}=\frac{1}{1-r^4}
\label{B5}
\end{equation}

The following equation is applied to Eq. (\ref{B5}).

\begin{equation}
(\mathrm{sin}(2\theta))^2=1-(\mathrm{cos}(2\theta))^2=1-r^4
\label{B6}
\end{equation}

Hence, the arc length $l$ becomes

\begin{equation}
l= \int_{r}^{1} \frac{1}{\sqrt{1-t^4}} \mathrm{d}t
\label{B7}
\end{equation}

\appendix
\def\thesection{Appendix}
\section{C}
\label{C}
The following function is differentiated.

\begin{equation}
\frac{\mathrm{d}}{\mathrm{d}l} \mathrm{arctan} (\mathrm{cleaf}_2(l))
= - \frac{ \sqrt{1- (\mathrm{cleaf}_2(l) )^4}  }{ 1+ (\mathrm{cleaf}_2(l) )^2  }
= - \sqrt{ \frac{ 1- (\mathrm{cleaf}_2(l) )^2  }{  1+ (\mathrm{cleaf}_2(l) )^2  }  }
= - \mathrm{sleaf}_2(l)  \label{C1}
\end{equation}

Integration of the abovementioned equation with respect to $l$  yields

\begin{equation}
\begin{split}
& \int_{0}^{l}  \mathrm{sleaf}_2(l) \mathrm{d}t
= \left[ -\mathrm{arctan} (\mathrm{cleaf}_2(l)) \right]^l_0
=-\mathrm{arctan} (\mathrm{cleaf}_2(l)) + \mathrm{arctan} (\mathrm{cleaf}_2(0)) \\
&=-\mathrm{arctan} (\mathrm{cleaf}_2(l)) + \mathrm{arctan} (1)
=-\mathrm{arctan} (\mathrm{cleaf}_2(l)) + \frac{\pi}{4}.
\label{C2}
\end{split}
\end{equation}

Similarly, the following function is differentiated with respect to variable $l$.

\begin{equation}
\frac{\mathrm{d}}{\mathrm{d}l} \mathrm{arctan} (\mathrm{sleaf}_2(l))
=  \frac{ \sqrt{1- (\mathrm{sleaf}_2(l) )^4}  }{ 1+ (\mathrm{sleaf}_2(l) )^2  }
=  \sqrt{ \frac{ 1- (\mathrm{sleaf}_2(l) )^2  }{  1+ (\mathrm{sleaf}_2(l) )^2  }  }
= \mathrm{cleaf}_2(l)  \label{C3}
\end{equation}

The following equation is obtained by integrating Eq. (\ref{C3}) with respect to $l$.

\begin{equation}
\begin{split}
& \int_{0}^{l}  \mathrm{cleaf}_2(l) \mathrm{d}t
= \left[ \mathrm{arctan} (\mathrm{sleaf}_2(l)) \right]^l_0
= \mathrm{arctan} (\mathrm{sleaf}_2(l)) - \mathrm{arctan} (\mathrm{sleaf}_2(0)) \\
&=\mathrm{arctan} (\mathrm{sleaf}_2(l)) - \mathrm{arctan} (0)
=\mathrm{arctan} (\mathrm{sleaf}_2(l))
\label{C4}
\end{split}
\end{equation}

\appendix
\def\thesection{Appendix}
\section{D}
\label{D}
Equation (\ref{4.1}) of the Cartesian coordinate system is transformed to the following equation in the polar coordinates \cite{Kaz_sl} \cite{Cox}.

\begin{equation}
\overline{r}^2=\mathrm{sin}(2 \overline{\theta}) \label{D1}
\end{equation}

By differentiating Eq. (\ref{D1}) with respect to the variable $\overline{\theta}$,

\begin{equation}
2\overline{r} \frac{\mathrm{d} \overline{r}}{\mathrm{d} \overline{\theta}} =2 \mathrm{cos}(2\overline{\theta})
\label{D4}
\end{equation}

Then, applying Eq (\ref{D4}),

\begin{equation}
1+\overline{r}^2 \left( \frac{ \mathrm{d} \overline{\theta}}{ \mathrm{d} \overline{r}} \right)^2=1+\overline{r}^2 \left( \frac{ \overline{r}}{\mathrm{cos}(2 \overline{\theta})} \right)^2 =1+\frac{\overline{r}^4}{1-\overline{r}^4}=\frac{1}{1-\overline{r}^4}
\label{D5}
\end{equation}

Hence, the arc length $\overline{l}$ becomes

\begin{equation}
\overline{l}= \int_{0}^{\overline{r}} \frac{1}{\sqrt{1-t^4}} \mathrm{d}t
\label{D7}
\end{equation}

\end{document}